%% file: IMEXRP.tex
\documentclass[final,preprint,3p,sort&compress]{elsarticle}
\journal{{\tt arXiv.org}}

\usepackage[dvips]{epsfig}
\usepackage{graphicx} 
\usepackage{pdfpages}
\usepackage{latexsym}
\usepackage{verbatim}
\usepackage{amsmath,amsthm}
\usepackage{amssymb}
\usepackage{bbm}
\usepackage{fonts}
\usepackage{enumerate}
\usepackage{placeins}
\usepackage{standalone}
\usepackage{paralist}
\usepackage{subcaption}
\definecolor{greenyellow}   {cmyk}{0.15, 0   , 0.69, 0   }
\definecolor{yellow}        {cmyk}{0   , 0   , 1   , 0   }
\definecolor{goldenrod}     {cmyk}{0   , 0.10, 0.84, 0   }
\definecolor{dandelion}     {cmyk}{0   , 0.29, 0.84, 0   }
\definecolor{apricot}       {cmyk}{0   , 0.32, 0.52, 0   }
\definecolor{peach}         {cmyk}{0   , 0.50, 0.70, 0   }
\definecolor{melon}         {cmyk}{0   , 0.46, 0.50, 0   }
\definecolor{yelloworange}  {cmyk}{0   , 0.42, 1   , 0   }
\definecolor{orange}        {cmyk}{0   , 0.61, 0.87, 0   }
\definecolor{burntorange}   {cmyk}{0   , 0.51, 1   , 0   }
\definecolor{bittersweet}   {cmyk}{0   , 0.75, 1   , 0.24}
\definecolor{redorange}     {cmyk}{0   , 0.77, 0.87, 0   }
\definecolor{mahogany}      {cmyk}{0   , 0.85, 0.87, 0.35}
\definecolor{maroon}        {cmyk}{0   , 0.87, 0.68, 0.32}
\definecolor{brickred}      {cmyk}{0   , 0.89, 0.94, 0.28}
\definecolor{red}           {cmyk}{0   , 1   , 1   , 0   }
\definecolor{orangered}     {cmyk}{0   , 1   , 0.50, 0   }
\definecolor{rubinered}     {cmyk}{0   , 1   , 0.13, 0   }
\definecolor{wildstrawberry}{cmyk}{0   , 0.96, 0.39, 0   }
\definecolor{salmon}        {cmyk}{0   , 0.53, 0.38, 0   }
\definecolor{carnationpink} {cmyk}{0   , 0.63, 0   , 0   }
\definecolor{magenta}       {cmyk}{0   , 1   , 0   , 0   }
\definecolor{violetred}     {cmyk}{0   , 0.81, 0   , 0   }
\definecolor{rhodamine}     {cmyk}{0   , 0.82, 0   , 0   }
\definecolor{mulberry}      {cmyk}{0.34, 0.90, 0   , 0.02}
\definecolor{redviolet}     {cmyk}{0.07, 0.90, 0   , 0.34}
\definecolor{fuchsia}       {cmyk}{0.47, 0.91, 0   , 0.08}
\definecolor{lavender}      {cmyk}{0   , 0.48, 0   , 0   }
\definecolor{thistle}       {cmyk}{0.12, 0.59, 0   , 0   }
\definecolor{orchid}        {cmyk}{0.32, 0.64, 0   , 0   }
\definecolor{darkorchid}    {cmyk}{0.40, 0.80, 0.20, 0   }
\definecolor{purple}        {cmyk}{0.45, 0.86, 0   , 0   }
\definecolor{plum}          {cmyk}{0.50, 1   , 0   , 0   }
\definecolor{violet}        {cmyk}{0.79, 0.88, 0   , 0   }
\definecolor{royalpurple}   {cmyk}{0.75, 0.90, 0   , 0   }
\definecolor{blueviolet}    {cmyk}{0.86, 0.91, 0   , 0.04}
\definecolor{periwinkle}    {cmyk}{0.57, 0.55, 0   , 0   }
\definecolor{cadetblue}     {cmyk}{0.62, 0.57, 0.23, 0   }
\definecolor{cornflowerblue}{cmyk}{0.65, 0.13, 0   , 0   }
\definecolor{midnightblue}  {cmyk}{0.98, 0.13, 0   , 0.43}
\definecolor{navyblue}      {cmyk}{0.94, 0.54, 0   , 0   }
\definecolor{royalblue}     {cmyk}{1   , 0.50, 0   , 0   }
\definecolor{blue}          {cmyk}{1   , 1   , 0   , 0   }
\definecolor{cerulean}      {cmyk}{0.94, 0.11, 0   , 0   }
\definecolor{cyan}          {cmyk}{1   , 0   , 0   , 0   }
\definecolor{processblue}   {cmyk}{0.96, 0   , 0   , 0   }
\definecolor{skyblue}       {cmyk}{0.62, 0   , 0.12, 0   }
\definecolor{turquoise}     {cmyk}{0.85, 0   , 0.20, 0   }
\definecolor{tealblue}      {cmyk}{0.86, 0   , 0.34, 0.02}
\definecolor{aquamarine}    {cmyk}{0.82, 0   , 0.30, 0   }
\definecolor{bluegreen}     {cmyk}{0.85, 0   , 0.33, 0   }
\definecolor{emerald}       {cmyk}{1   , 0   , 0.50, 0   }
\definecolor{junglegreen}   {cmyk}{0.99, 0   , 0.52, 0   }
\definecolor{seagreen}      {cmyk}{0.69, 0   , 0.50, 0   }
\definecolor{green}         {cmyk}{1   , 0   , 1   , 0   }
\definecolor{forestgreen}   {cmyk}{0.91, 0   , 0.88, 0.12}
\definecolor{pinegreen}     {cmyk}{0.92, 0   , 0.59, 0.25}
\definecolor{limegreen}     {cmyk}{0.50, 0   , 1   , 0   }
\definecolor{yellowgreen}   {cmyk}{0.44, 0   , 0.74, 0   }
\definecolor{springgreen}   {cmyk}{0.26, 0   , 0.76, 0   }
\definecolor{olivegreen}    {cmyk}{0.64, 0   , 0.95, 0.40}
\definecolor{rawsienna}     {cmyk}{0   , 0.72, 1   , 0.45}
\definecolor{sepia}         {cmyk}{0   , 0.83, 1   , 0.70}
\definecolor{brown}         {cmyk}{0   , 0.81, 1   , 0.60}
\definecolor{tan}           {cmyk}{0.14, 0.42, 0.56, 0   }
\definecolor{gray}          {cmyk}{0   , 0   , 0   , 0.50}
\definecolor{black}         {cmyk}{0   , 0   , 0   , 1   }
\definecolor{white}         {cmyk}{0   , 0   , 0   , 0   } 

\usepackage{booktabs} % for \cmidrule in tables
\usepackage{pgfplots}
\pgfplotsset{compat=newest}

\newcommand{\externaltikz}[2]{\includegraphics{Externals/#1}}	
 \usepackage{relinput}
\newtheorem{theorem}{Theorem}[section]
\newtheorem{definition}[theorem]{Definition}
\newtheorem{remark}[theorem]{Remark}
\newtheorem{example}[theorem]{Example}
\newtheorem{assumption}[theorem]{Assumption}

\newtheorem{lemma}[theorem]{Lemma}
\newtheorem{corollary}[theorem]{Corollary}

\newcounter{tikzsubfigcounter}[figure]
\renewcommand{\thetikzsubfigcounter}{\the\numexpr\value{figure}+1\relax\alph{tikzsubfigcounter}}

\newcounter{tikzsubfigcounterinvisible}[figure]
\renewcommand{\thetikzsubfigcounterinvisible}{\the\numexpr\value{figure}+1\relax\alph{tikzsubfigcounterinvisible}}

\newcommand{\settikzlabel}[1]{ %
\refstepcounter{tikzsubfigcounterinvisible} \label{#1} 
}

\numberwithin{equation}{section}
\title{Implicit-explicit, realizability-preserving first-order scheme for moment models with Lipschitz-continuous source terms}

\author[fs]{Florian Schneider}
\address[fs]{Fachbereich Mathematik, TU Kaiserslautern, Erwin-Schr\"odinger-Str., 67663 Kaiserslautern, Germany, {\tt schneider@mathematik.uni-kl.de}}

\date{}

\input{usermacros}

\begin{document}

\begin{abstract}
We derive an implicit-explicit (IMEX), realizability-preserving first-order scheme for moment models with Lipschitz-continuous source terms. In contrast to the fully-explicit schemes in \cite{Schneider2015a,Schneider2015b} the time step does not depend on the physical parameters, removing the stiffness from the system. Furthermore, a wider class of collision operators (e.g. the Laplace-Beltrami operator) can be used. The derived scheme is applied to minimum-entropy models.
\end{abstract}
\maketitle

\noindent

% {\bf Key words.}

\input{Sections/introduction}
\input{Sections/modelling}
\input{Sections/realizabilityred}
\input{Sections/scheme}
\input{Sections/results}

\input{Sections/outlook}

% Bibliography
%%%%%%%%%%%%%%
\bibliographystyle{siam}
\bibliography{bibliography}

\end{document}

%% file: usermacros.tex
\newlength{\figureheight}
\newlength{\figurewidth}
\usetikzlibrary{arrows,shapes,positioning}
\usetikzlibrary{decorations.markings}
\tikzstyle arrowstyle=[scale=1]
\tikzstyle directed=[postaction={decorate,decoration={markings,
		mark=at position .65 with {\arrow[arrowstyle]{stealth}}}}]
\tikzstyle reverse directed=[postaction={decorate,decoration={markings,
		mark=at position .65 with {\arrowreversed[arrowstyle]{stealth};}}}]

\newcommand{\secref}[1]{Section~\ref{#1}}
\newcommand{\thmref}[1]{Theorem~\ref{#1}}

\newcommand{\lemref}[1]{Lemma~\ref{#1}}

\renewcommand{\corref}[1]{Corollary~\ref{#1}}
\newcommand{\figref}[1]{Figure~\ref{#1}}
\newcommand{\tabref}[1]{Table~\ref{#1}}

\newcommand{\assref}[1]{Assumption~\ref{#1}}

\newcommand{\abs}[1]{\ensuremath{\left| #1 \right|}}
\newcommand{\norm}[2]{\ensuremath{\left\Vert #1 \right\Vert_{#2}}}

\newcommand{\N}{\mathbb{N}} 

\newcommand{\R}{\mathbb{R}}
\newcommand{\Rpos}{\R_{\geq 0}}
\newcommand{\scattering}{\ensuremath{\sigma_s}}
\newcommand{\absorption}{\ensuremath{\sigma_a}}

\newcommand{\source}{\ensuremath{Q}}

\newcommand{\distribution}[1][ ]{\ensuremath{\psi_{#1}}}
\newcommand{\distributiontzero}{\ensuremath{\distribution[\timevar=0]}}
\newcommand{\distributionboundary}{\ensuremath{\distribution[b]}}
\newcommand{\distributionvacuum}{\ensuremath{\distribution[\text{vac}]}}
\newcommand{\ansatz}[1][ ]{\ensuremath{\hat{\psi}_{#1}}}

\newcommand{\momentorder}{\ensuremath{N}}
\newcommand{\momentnumber}{\ensuremath{n}}
\newcommand{\basis}[1][ ]{{\ensuremath{\bb_{#1}}}} %Basis
\newcommand{\basisind}{\ensuremath{i}} %Basis
\newcommand{\basiscomp}[1][\basisind]{\ensuremath{b_{#1}}} %Basis
\newcommand{\fmbasis}[1][\momentorder]{\basis[#1]} %Basis
\newcommand{\moments}[1][ ]{\ensuremath{\bu_{#1}}} %moment vector
\newcommand{\momentcomp}[1]{\ensuremath{u_{#1}}} %moment vector

\newcommand{\normalizedmoments}[1][ ]{\ensuremath{\bsphi_{#1}}} %moment vector
\newcommand{\normalizedmomentcomp}[1]{\ensuremath{\phi_{#1}}} %moment vector
\newcommand{\normalizedisotropicmoment}{\normalizedmoments[\text{iso}]}
\newcommand{\multipliers}[1][ ]{\ensuremath{\bsalpha_{#1}}} %moment vector
\newcommand{\multiplierscomp}[1]{\ensuremath{\alpha_{#1}}} %moment vector
\newcommand{\SC}{\ensuremath{\Omega}} %moment vector
\newcommand{\SCheight}{\ensuremath{\mu}} %moment vector
\newcommand{\LaplaceBeltrami}{\ensuremath{\Delta_\SC}} %moment vector
\newcommand{\Domain}{\ensuremath{X}} %moment vector
\newcommand{\spatialVariable}{\ensuremath{\bx}} %moment vector
\newcommand{\timeint}{\ensuremath{T}} %time interval
\newcommand{\tf}{\ensuremath{t_f}} %final time
\newcommand{\timevar}{\ensuremath{t}} %final time

\newcommand{\ints}[1]{\ensuremath{\left<#1\right>}}

\newcommand{\collisionop}{\ensuremath{\cC}}
\newcommand{\collisionopu}{\ensuremath{\bC}}
\newcommand{\collision}[1]{\ensuremath{\collisionop\left(#1\right)}}
\newcommand{\collisionu}[1]{\ensuremath{\collisionopu\left(#1\right)}}

\newcommand{\collisionkernel}{\ensuremath{K}}
\newcommand{\collisionrealizablepart}[1][ ]{\ensuremath{\widetilde{\moments[#1]}}}

\newcommand{\indicator}[1]{\ensuremath{\mathbbm{1}_{#1}}}

\newcommand{\Lp}[1]{\ensuremath{L_{#1}}}
\newcommand{\RD}[2]{\ensuremath{\mathcal{R}_{#1}^{#2}}}
\newcommand{\RDone}[1]{\left. \RD{#1}{} \right|_{\momentcomp{0} = 1}}
\newcommand{\dRDone}[1]{\left. \partial\RD{#1}{} \right|_{\momentcomp{0} = 1}}

\newcommand{\distRD}[1]{\ensuremath{d\left(#1,\dRDone{\basis}\right)}}
\newcommand{\distRDrel}[1]{\ensuremath{d_r\left(#1,\dRDone{\basis}\right)}}

\newcommand{\AnsatzSpace}{\ensuremath{\cA}}

\newcommand{\PN}[1][\momentorder]{\ensuremath{\text{P}_{#1}}}
\newcommand{\MN}[1][\momentorder]{\ensuremath{\text{M}_{#1}}}

\newcommand{\Flux}{\ensuremath{\bF}}
\newcommand{\Source}{\ensuremath{\bs}}

\DeclareMathOperator*{\argmin}{argmin}
\newcommand{\ld}[1]{\ensuremath{{#1}_*}} %Legendre dual
\newcommand{\entropy}{\ensuremath{\eta}} %Entropy function
\newcommand{\entropyFunctional}{\ensuremath{\mathcal{H}}} %Entropy functional
\def\quand{\quad \mbox{and} \quad}

\newcommand{\x}{\ensuremath{x}}

\newcommand{\z}{\x}
\newcommand{\dx}{\partial_{\x}}

\newcommand{\dz}{\partial_{\z}}
\newcommand{\dt}{\partial_\timevar}

\newcommand{\LaplaceBeltramiProjection}{\ensuremath{\Delta_\SCheight}} %moment vector

\newcommand{\SobolevToyproblem}{\ensuremath{\mathcal{V}}}
\newcommand{\SobolevTestfunction}{\ensuremath{v}}

\newcommand{\dtstepsize}{\ensuremath{\Delta \timevar}}

\newcommand{\opttol}{\ensuremath{\tau}}

\newcommand{\zL}{\ensuremath{\z_{L}}}
\newcommand{\zR}{\ensuremath{\z_{R}}}
\newcommand{\dzstepsize}{\ensuremath{\Delta\z}}
\newcommand{\ncells}{\ensuremath{n_{\z}}}

\newcommand{\timeind}{\ensuremath{\kappa}}

\newcommand{\cellind}{\ensuremath{j}}

\newcommand{\cell}[1]{\ensuremath{I_{#1}}}

\newcommand{\cellmean}[2][\cellind]{\ensuremath{\overline{#2}_{#1}}}

\newcommand{\monotoneFluxSymbol}{\ensuremath{h}}
\newcommand{\monotoneFlux}[2]{\ensuremath{\monotoneFluxSymbol\left(#1,#2\right)}}

\newcommand{\distributionpv}[1]{\ensuremath{\distribution[#1]}}

\newcommand{\momentsprojected}{\ensuremath{\moments[h]}}
\newcommand{\momentcompprojected}[1]{\ensuremath{\momentcomp{h,#1}}}

\newcommand{\momentscellmeantime}[2]{\ensuremath{\cellmean[#1]{\moments}^{\left(#2\right)}}}
\newcommand{\momentslocal}[1]{\moments[#1]}

\newcommand{\numericalFlux}{\ensuremath{\widehat{\bF}}}

\newcommand{\LpError}[2][h]{\ensuremath{E_{#1}^{#2}}}

\newcommand{\MFSconstK}{\ensuremath{K}}
\newcommand{\MFSconsta}{\ensuremath{b}}

\def\analyticalsolution{\distribution[a]}
\def\analyticalmoments{\moments[a]}
\newcommand{\analyticalmomentcomp}[1]{\ensuremath{\momentcomp{a,#1}}}
\newcommand{\MFSconstorder}{\ensuremath{\nu}}
\newcommand{\MFSp}{\ensuremath{p}}

\newcommand{\density}{\ensuremath{\rho}}

\pgfplotscreateplotcyclelist{color parula}{% 
	royalblue!70,every mark/.append style={solid,line width = 0pt,fill=royalblue!60!black},mark=ball\\%
	goldenrod!50!yelloworange,every mark/.append style={solid,fill=goldenrod!30!black},mark=*\\%
	green,every mark/.append style={black,fill=yellowgreen},mark=diamond*\\%
	bluegreen,every mark/.append style={fill=black},mark=triangle*\\%	
	royalblue!70,densely dashed,every mark/.append style={solid,line width = 0pt,fill=royalblue!60!black},mark=ball\\%
	goldenrod!50!yelloworange,densely dashed,every mark/.append style={solid,black,fill=goldenrod},mark=diamond*\\%
	yellowgreen,densely dashed,every mark/.append style={solid,fill=yellowgreen!60!royalblue},mark=square*, mark size= 1pt\\%
	bluegreen,densely dashed,mark size = 1.5pt,every mark/.append style={solid,fill=white},mark=otimes*\\%,densely dashed	
	royalblue,densely dashed,mark size = 1.5pt,every mark/.append style={solid,fill=royalblue!30!black},mark=pentagon*\\%,densely dashed
	goldenrod!40!yelloworange,every mark/.append style={solid,fill=goldenrod!30!black},mark=|\\%		
	}

%% file: Sections/introduction.tex
\section{Introduction}
In recent years many approaches have been considered for the solution of
time-dependent linear kinetic transport equations, which arise for example
in electron radiation therapy or radiative heat transfer problems.
Many of the most popular methods are moment methods, also known as moment
closures because they are distinguished by how they close the truncated
system of exact moment equations.
Moments are defined through angular averages against basis functions to
produce spectral approximations in the angle variable.
A typical family of moment models are the so-called $\PN$-methods
\cite{Lewis-Miller-1984,Gel61} which are pure spectral methods.
However, many high-order moment methods, including $\PN$, do not take
into account that the original kinetic density to be approximated must be
non-negative.
The moment vectors produced by such models are therefore often not realizable,
that is, there is no associated non-negative kinetic distribution consistent with
the moment vector, and thus the solutions can contain non-physical
artefacts such as negative local particle densities \cite{Bru02}.

The family of minimum-entropy models, colloquially known as $\MN$ models
or entropy-based moment closures, solve this problem (for certain physically
relevant entropies) by specifying the closure using a non-negative density
reconstructed from the moments.
The $\MN$ models are the only models which additionally are hyperbolic and
dissipate entropy \cite{Lev96}.
The cost of all these properties is that the reconstruction of this density
involves solving an optimization problem at every point on the space-time
mesh \cite{AllHau12,Alldredge2014}.
These reconstructions, however, can be parallelized, and so the recent emphasis
on algorithms that can take advantage of massively parallel computing
environments has led to renewed interest in the computation of $\MN$ solutions
both for linear and nonlinear kinetic equations
\cite{DubFeu99,Hauck2010,Lam2014,AllHau12,Garrett2014,
McDonald2012}.

The key challenge for a numerical scheme is that, if not treated correctly, the numerical solution can leave the set of realizable moments \cite{Olbrant2012}, outside of which the defining optimization problem has no solution.

Discontinuous-Galerkin methods can handle this problem using a realizability
limiter directly on the moment vectors themselves
\cite{Zhang2010,Olbrant2012,Schneider2015a}. At this level
realizability conditions are in general quite complicated
and also not well-understood for two- or three-dimensional problems for
moment models of order higher than two.
Realizability limiting for kinetic schemes \cite{Hauck2010,Schneider2015b}, however, is much easier because
at the level of the kinetic density, realizability corresponds simply to
non-negativity.

One big drawback of explicit schemes is that the time step depends on the physical parameters (absorption and scattering properties of the material), resulting in stiff systems, which can be avoided using an implicit discretization.
On the other hand, the hyperbolic flux, which is non-linear and usually expensive to calculate, is typically non-stiff. An implicit discretization is therefore undesired. To overcome this we derive a realizability-preserving, first-order kinetic scheme with implicit-explicit (IMEX) time stepping, treating stiff and non-stiff problems separately.

The paper is organized as follows. A brief overview of the method of moment, the minimum-entropy approach and realizability is given in \secref{sec:Models}. Then, the reduced (space-homogeneous) moment system (which will be treated implicitly in the scheme) is investigated and the realizability-preserving property of this implicit discretization is shown in \secref{sec:RealizabilityReduced}. This is concluded by the description of the full scheme and the proof that it is realizability-preserving in \secref{sec:RPFO}. The scheme is then tested in a manufactured solution and a benchmark test in \secref{sec:NumExp}. Finally, conclusions and an outlook on future work is given in \secref{sec:Conclusions}.

%% file: Sections/modelling.tex
\section{Models}
\label{sec:Models}
In slab geometry, the transport equation under consideration has the form 
\begin{align}
\label{eq:TransportEquation1D}
\dt\distribution+\SCheight\dz\distribution + \absorption\distribution = \scattering\collision{\distribution}+\source, \qquad \timevar\in\timeint,\z\in\Domain,\SCheight\in[-1,1].
\end{align}
The physical parameters are the absorption and scattering coefficient $\absorption,\scattering:\timeint\times\Domain\to\Rpos$, respectively, and the emitting source $\source:\timeint\times\Domain\times[-1,1]\to\Rpos$. Furthermore, $\SCheight\in[-1,1]$, and $\distribution = \distribution(\timevar,\z,\SCheight)$. 

\begin{assumption}
\label{ass:CollisionOperator}
The operator $\collisionop$ is assumed to have the following properties.
\begin{enumerate}
\begin{subequations}
\label{eq:CollisionProperty}
\item Mass conservation
\begin{align}
\label{eq:CollisionPropertyMass}
\int\limits_{-1}^1\collision{\distribution}~d\SCheight=0.
\end{align}
\item Local entropy dissipation
\begin{align}
\label{eq:CollisionPropertyLocalDissipation}
\int\limits_{-1}^1\entropy'(\distribution)\collision{\distribution}~d\SCheight\leq 0,
\end{align}
where $\entropy$ denotes a strictly convex entropy (compare \secref{sec:MinimumEntropy}).
\item The reduced (space-homogeneous) system 
\begin{align}
\label{eq:CollisionDistributionEquation}
\dt \distribution = \collision{\distribution},
\end{align}
admits a non-negative solution $\distribution\geq 0$ for all $\timevar\geq 0$ and initial conditions $\distribution(0,\SCheight)\geq 0$.
\item For every $\dtstepsize\geq0$, the following implication holds
\begin{align}
\label{eq:CollisionDistributionPositivity}
\distribution(\timevar,\SCheight)-\dtstepsize\collision{\distribution(\timevar,\SCheight)}\geq 0~\Rightarrow~ \distribution(\timevar,\SCheight)\geq 0 \quad\text{for all } \timevar\geq 0,~\SCheight\in[-1,1].
\end{align}
\end{subequations}
\end{enumerate}
\end{assumption}
The first two assumptions are from \cite{Levermore1996}, requiring that the operator is physically meaningful. The other assumptions are necessary for some of our proofs in the following\footnote{To be completely correct, the assumptions have to be formulated in a weak sense. However, all steps below can be performed similarly but with a greater notational effort.}. 
One example for such a collision operator is given by the Laplace-Beltrami operator 
\begin{align}
\label{eq:LaplaceBeltrami}
\collision{\distribution} = \frac12 \LaplaceBeltramiProjection \distribution = \frac12 \cfrac{d}{d\SCheight}\left(\left(1-\SCheight^2\right)\cfrac{d\distribution}{d\SCheight}\right).
\end{align}
This operator appears, for example, as the result of an asymptotic analysis of the Boltzmann equation under the assumption of small energy loss and deflection, and forward-peaked scattering in the context of electron transport \cite{Frank07,Pom92,HenIzaSie06}.

Another typical choice is the linear integral collision operator 
\begin{equation}
 \collision{\distribution} =  \int\limits_{-1}^1 \collisionkernel(\SCheight, \SCheight^\prime)
  \distribution(\timevar, \spatialVariable, \SCheight^\prime)~d\SCheight^\prime 
  - \int\limits_{-1}^1 \collisionkernel(\SCheight^\prime, \SCheight) \distribution(\timevar, \spatialVariable, \SCheight)~d\SCheight^\prime.
\label{eq:collisionOperatorR}
\end{equation}
The collision kernel $\collisionkernel$ is assumed to be strictly positive, symmetric (i.e. $\collisionkernel(\SCheight,\SCheight')=\collisionkernel(\SCheight',\SCheight)$) and normalized to 
$\int\limits_{-1}^1 \collisionkernel(\SCheight^\prime, \SCheight)~d\SCheight^\prime=1$.  A typical example is the BGK-type
\emph{isotropic-scattering} operator, where $\collisionkernel(\SCheight, \SCheight^\prime) \equiv \frac{1}{2}$.

The transport equation \eqref{eq:TransportEquation1D} is supplemented by initial and boundary conditions:
\begin{subequations}
\begin{align}
\distribution(0,\z,\SCheight) &= \distributiontzero(\z,\SCheight) &&\text{for } \z\in\Domain = (\zL,\zR), \SCheight\in[-1,1], \label{eq:TransportEquation1DIC}\\
\distribution(\timevar,\zL,\SCheight) &= \distributionboundary(\timevar,\zL,\SCheight) &&\text{for } \timevar\in\timeint, \SCheight>0,  \label{eq:TransportEquation1DBCa}\\
\distribution(\timevar,\zR,\SCheight) &= \distributionboundary(\timevar,\zR,\SCheight) &&\text{for } \timevar\in\timeint, \SCheight<0. \label{eq:TransportEquation1DBCb}
\end{align}
\end{subequations}

\subsection{The method of moments}
In general, solving equation \eqref{eq:TransportEquation1D} is very expensive in two and three dimensions due to the high dimensionality of the state space. 

For this reason it is convenient to use some type of spectral or Galerkin method to transform the high-dimensional equation into a system of lower-dimensional equations. Typically, one chooses to reduce the dimensionality by representing the angular dependence of $\distribution$ in terms of some basis $\basis$.
\begin{definition}
The vector of functions $\basis:[-1,1]\to\R^{\momentnumber}$ consisting of $\momentnumber$ basis functions $\basiscomp[\basisind]$, $\basisind=0,\ldots\momentnumber-1$ of maximal \emph{order} $\momentorder$ is called an \emph{angular basis}. Analogously, the symbol $\basis[\momentorder]$ can be used if the knowledge of $\momentorder$ is explicitly necessary.

The so-called \emph{moments} of a given distribution function $\distribution$ with respect to $\basis$ are then defined by
\begin{align}
\label{eq:moments}
\moments =\ints{{\basis}\distribution} = \left(\momentcomp{0},\ldots,\momentcomp{\momentnumber-1}\right)^T,
\end{align}
where the integration $\ints{\cdot} = \int\limits_{-1}^1\cdot~d\SCheight$ is performed componentwise.\\

Assuming for simplicity $\basiscomp[0]\equiv 1$, the quantity $\momentcomp{0} = \ints{\basiscomp[0]\distribution}=\ints{\distribution}$ is called \emph{local particle density}. 
Furthermore, \emph{normalized moments} $\normalizedmoments = \left(\normalizedmomentcomp{1},\ldots,\normalizedmomentcomp{\momentnumber-1}\right)\in\R^{\momentnumber}$ are defined as 
\begin{align}
\label{eq:NormalizedMoments}
\normalizedmomentcomp{\basisind} = \cfrac{\momentcomp{\basisind}}{\momentcomp{0}}~, \qquad \basisind=1,\ldots\momentnumber-1.
\end{align}
\end{definition}
To obtain a set of equations for $\moments$, \eqref{eq:TransportEquation1D} has to be multiplied through by $\basis$ and integrated over $[-1,1]$, giving
\begin{align*}
\ints{\basis\dt\distribution}+\ints{\basis\dz\SCheight\distribution} + \ints{\basis\absorption\distribution} = \scattering\ints{\basis\collision{\distribution}}+\ints{\basis\source}.
\end{align*}
Collecting known terms, and interchanging integrals and differentiation where possible, the moment system has the form
\begin{align}
\label{eq:MomentSystemUnclosed1D}
\dt\moments+\dz\ints{\SCheight \basis\ansatz[\moments]} + \absorption\moments = \scattering\ints{\basis\collision{\ansatz[\moments]}}+\ints{\basis\source}.
\end{align}

The solution of \eqref{eq:MomentSystemUnclosed1D} is equivalent to the one of \eqref{eq:TransportEquation1D} if $\basis$ is a basis of $\Lp{2}([-1,1],\R)$. 

Since it is impractical to work with an infinite-dimensional system, only a finite number of $\momentnumber<\infty$ basis functions $\basis$ of order $\momentorder$ can be considered. Unfortunately, there always exists an index $\basisind\in\{0,\dots,\momentnumber-1\}$ such that the components of $\basiscomp\cdot\SCheight$ are not in the linear span of $\basis$. Therefore, the flux term cannot be expressed in terms of $\moments$ without additional information. Furthermore, the same might be true for the projection of the scattering operator onto the moment-space given by $\ints{\basis\collision{\distribution}}$. This is the so-called \emph{closure problem}. One usually prescribes some \emph{ansatz} distribution $\ansatz[\moments](\timevar,\spatialVariable,\SCheight):=\ansatz(\moments(\timevar,\spatialVariable),\basis(\SCheight))$ to calculate the unknown quantities in \eqref{eq:MomentSystemUnclosed1D}. Note that the dependence on the angular basis in the short-hand notation $\ansatz[\moments]$ is neglected for notational simplicity.

Finally, we write \eqref{eq:MomentSystemUnclosed1D} in the form of a standard first-order hyperbolic system of equations:
\begin{align}
\label{eq:GeneralHyperbolicSystem}
\dt\moments + \dx\Flux(\moments) = \Source\left(\moments\right),
\end{align}
where $\Flux(\moments) = \ints{\SCheight \basis\ansatz[\moments]}$ and $\Source\left(\moments\right) = \scattering\ints{\basis\collision{\ansatz[\moments]}}+\ints{\basis\source}-\absorption\moments$.
\subsection{Minimum-entropy approach}
\label{sec:MinimumEntropy}
In this paper the ansatz density $\ansatz$ is reconstructed from the moments $\moments$ by minimizing the entropy-functional 
 \begin{align}
 \label{eq:entropyFunctional}
 \entropyFunctional(\distribution) = \ints{\entropy(\distribution)}
 \end{align}
 under the moment constraints
 \begin{align}
 \label{eq:MomentConstraints}
 \ints{\basis\distribution} = \moments.
 \end{align}
The kinetic entropy density $\entropy:\R\to\R$ is strictly convex and twice continuously differentiable and the minimum is simply taken over all functions $\distribution = \distribution(\SCheight)$ such that 
  $\entropyFunctional(\distribution)$ is well defined. The obtained ansatz $\ansatz = \ansatz[\moments]$, solving this constrained optimization problem, is given by
 \begin{equation}
  \ansatz[\moments] = \argmin\limits_{\distribution:\entropy(\distribution)\in\Lp{1}}\left\{\ints{\entropy(\distribution)}
  : \ints{\basis \distribution} = \moments \right\}.
 \label{eq:primal}
 \end{equation}
This problem, which must be solved over the space-time mesh, is typically solved through its strictly convex finite-dimensional dual,
 \begin{equation}
  \multipliers(\moments) := \argmin_{\tilde{\multipliers} \in \R^{\momentnumber}} \ints{\ld{\entropy}(\basis^T 
   \tilde{\multipliers})} - \moments^T \tilde{\multipliers},
 \label{eq:dual}
 \end{equation}
where $\ld{\entropy}$ is the Legendre dual of $\entropy$. The first-order necessary conditions for the multipliers $\multipliers(\moments)$ show that the solution to \eqref{eq:primal} has the form
 \begin{equation}
  \ansatz[\moments] = \ld{\entropy}' \left(\basis^T \multipliers(\moments) \right)
 \label{eq:psiME}
 \end{equation}
where $\ld{\entropy}'$ is the derivative of $\ld{\entropy}$.\\

This approach is called the \emph{minimum-entropy closure} \cite{Levermore1996}. The resulting model has many desirable properties: symmetric hyperbolicity, bounded eigenvalues of the directional flux Jacobian and the direct existence of an entropy-entropy flux pair (compare \cite{Levermore1996,Schneider2016}).\\

The kinetic entropy density $\entropy$ can be chosen according to the 
physics being modelled.
As in \cite{Levermore1996,Hauck2010}, Maxwell-Boltzmann entropy%
 \begin{align}
 \label{eq:EntropyM}
  \entropy(\distribution) = \distribution \log(\distribution) - \distribution
 \end{align}
is used, thus $\ld{\entropy}(p) = \ld{\entropy}'(p) = \exp(p)$. This entropy is used for non-interacting particles as in an ideal gas.

We use the modification of the adaptive-basis optimization routine \cite{Alldredge2014} as proposed in \cite{Schneider2015b} to solve \eqref{eq:dual}.\\

Substituting $\distribution$ in \eqref{eq:MomentSystemUnclosed1D} with $\ansatz[\moments]$ yields a closed system of equations for $\moments$:
\begin{align}
\label{eq:MomentSystemClosed}
\dt\moments+\dz\ints{\SCheight \basis\ansatz[\moments]} + \absorption\moments = \scattering\ints{\basis\collision{\ansatz[\moments]}}+\ints{\basis\source}.
\end{align}

In this paper, the full-moment basis $\basis = \left(1,\SCheight,\ldots,\SCheight^\momentorder\right)$ will be used. Nevertheless, the scheme can be transferred directly to other bases like the \emph{half-moment monomial basis} ($\basiscomp = \indicator{[-1,0]}\SCheight^\basisind$ or $\basiscomp = \indicator{[0,1]}\SCheight^\basisind$) \cite{DubKla02,DubFraKlaTho03,Ritter2016} or the mixed-moment basis ($\basis = \left(1,\SCheight \indicator{[0,1]},\ldots,\SCheight^\momentorder \indicator{[0,1]},\SCheight \indicator{[-1,0]},\ldots,\SCheight^\momentorder \indicator{[-1,0]}\right)$) \cite{Schneider2014,Frank07,Schneider2015c}. Similarly, the results are not restricted to the minimum-entropy approach but can be transferred to other realizable closures like Kershaw \cite{Ker76,Schneider2016a,Schneider2015} or the quadrature method of moments \cite{Yuan2012,Fox2008,Fox2009,Vikas2013}.

\subsection{Realizability}

Since the underlying kinetic density to be approximated is
non-negative, a 
moment vector only makes sense physically if it can be associated with a 
non-negative distribution function. In this case the moment vector is called 
\emph{realizable}.

\begin{definition}
\label{def:RealizableSet}
The \emph{realizable set} $\RD{\basis}{}$ is 
$$
\RD{\basis}{} = \left\{\moments~:~\exists \distribution(\SCheight)\ge 0,\, \momentcomp{0} = \ints{\distribution} > 0,
 \text{ such that } \moments =\ints{\basis\distribution} \right\}.
$$
If $\moments\in\RD{\basis}{}$, then $\moments$ is called \emph{realizable}.
Any $\distribution$ such that $\moments =\ints{\basis \distribution}$ is called a \emph{representing 
density}.
\end{definition}
\begin{remark}
\mbox{ }
\begin{enumerate}[(a)]
\item The realizable set is a convex cone, and
\item Representing densities are not necessarily unique.
\end{enumerate}
\end{remark}

Additionally, since the entropy ansatz has the form \eqref{eq:psiME}, in the 
Maxwell-Boltzmann case, the optimization problem \eqref{eq:primal} only has a 
solution if the moment vector lies in the ansatz space
$$
 \AnsatzSpace := \left\{\ints{\basis \ansatz[\moments]}\stackrel{\eqref{eq:psiME}}{=} \ints{\basis \ld{\entropy}'\left(\basis^T\multipliers\right) }
  : \multipliers \in \R^{\momentnumber}  \right\}.
$$
In the case of a bounded angular domain, the ansatz 
space $\AnsatzSpace$ is equal to the set of realizable moment vectors 
\cite{Jun00}. Therefore, it is sufficient to focus on realizable moments only.

The definition of the realizable set is not constructive, making it hard to check if a moment vector is realizable or not. There are several works about concrete representations of the realizable set for different bases, e.g. \cite{Curto1991,Ker76,Schneider2014,Schneider2015a,Curto1996,Schneider2015c}. 

Exemplarily, the full-moment realizable set of order $\momentorder=2$ is given by \cite{Curto1991}
\begin{align}
\label{eq:M2Realizability2}
\RD{\basis}{} = \left\{\moments\in\R^3~:~\momentcomp{0}\geq\abs{\momentcomp{1}},~\momentcomp{0}\momentcomp{2}\geq\momentcomp{1}^2\right\}.
\end{align}

Fortunately, since we only use a first-order scheme, no information about the realizable set (except its convexity) is needed in the following. Note that this might no longer be true when higher-order schemes (in space and time) are used, see e.g. \cite{Schneider2015b,Zhang2010,Zhang2011a}.

%% file: Sections/realizabilityred.tex
\section{Realizability of the reduced equation}
\label{sec:RealizabilityReduced}
Before treating the space-dependent transport equation \eqref{eq:TransportEquation1D}, we want to investigate \eqref{eq:CollisionDistributionEquation} in more detail. The following example shows why explicit schemes for the Laplace-Beltrami fail.

\begin{example}
\begin{align}
\label{eq:LaplaceBeltramiMomentsEquation}
\dt \moments = \ints{\basis\LaplaceBeltramiProjection\ansatz[\moments]},
\end{align}
where the ansatz $\ansatz[\moments]$ can be chosen accordingly as \eqref{eq:psiME}, if necessary.\\

It is possible to show that \eqref{eq:CollisionDistributionEquation} has a solution in $\Lp{2}([-1,1],\Rpos)$ for every $\timevar\geq 0$ \cite{Kuo2006,risken1996fokker,hsu2002stochastic,hackenbroch1994stochastische}. Therefore, it is possible to expand $\distribution$ in $\SCheight$ in terms of the Legendre polynomials $P_\basisind$, which form an orthogonal basis of $\Lp{2}$ and are eigenfunctions of $\LaplaceBeltramiProjection$. Then, \eqref{eq:CollisionDistributionEquation} transforms to
\begin{align*}
\sum\limits_{\basisind=0}^{\infty} \left(\dt\multiplierscomp{\basisind}+\basisind\left(\basisind+1\right)\multiplierscomp{\basisind}\right)c_\basisind P_\basisind = 0,
\end{align*}
where the coefficients $c_\basisind$ are normalization constants. This equation can be stated equivalently as an infinite, decoupled system of ordinary differential equations
\begin{align*}
\dt\multiplierscomp{\basisind} = -\basisind\left(\basisind+1\right)\multiplierscomp{\basisind},\qquad\quad \basisind \in \N_{\geq 0}
\end{align*}
with solution
\begin{align*}
\multiplierscomp{\basisind}(\timevar) = e^{-\basisind\left(\basisind+1\right)\timevar}\multiplierscomp{\basisind}(0),
\end{align*}
where $\multiplierscomp{\basisind}(0)$ are the Fourier coefficients of $\distribution(0,\SCheight)$. For $\timevar\to\infty$ it obviously holds that 
\begin{align*}
\lim\limits_{\timevar\to\infty}\multiplierscomp{\basisind}(\timevar) = 0,&& \basisind = 1,\ldots,\infty
\end{align*}
which means that $\distribution(\timevar,\SCheight)\stackrel{\timevar\to\infty}{\longrightarrow} \multiplierscomp{0}(0)$. This implies that for every initial condition for \eqref{eq:CollisionDistributionEquation} a stationary solution is attained and that it is isotropic. This is not very surprising since the constants are in the kernel of $\LaplaceBeltramiProjection$.\\

The corresponding second-order, full-moment vector field 
\begin{align}
\label{eq:LaplaceBeltramiMomentsEquationFM}
\ints{\basis[2]\LaplaceBeltrami\ansatz} = (0,-2\momentcomp{1},-6\momentcomp{2}+2\momentcomp{0})^T
\end{align}
is plotted in normalized moments in \figref{fig:M2LaplaceBeltramiVectorField}. 

\begin{figure}[h]
\centering
\externaltikz{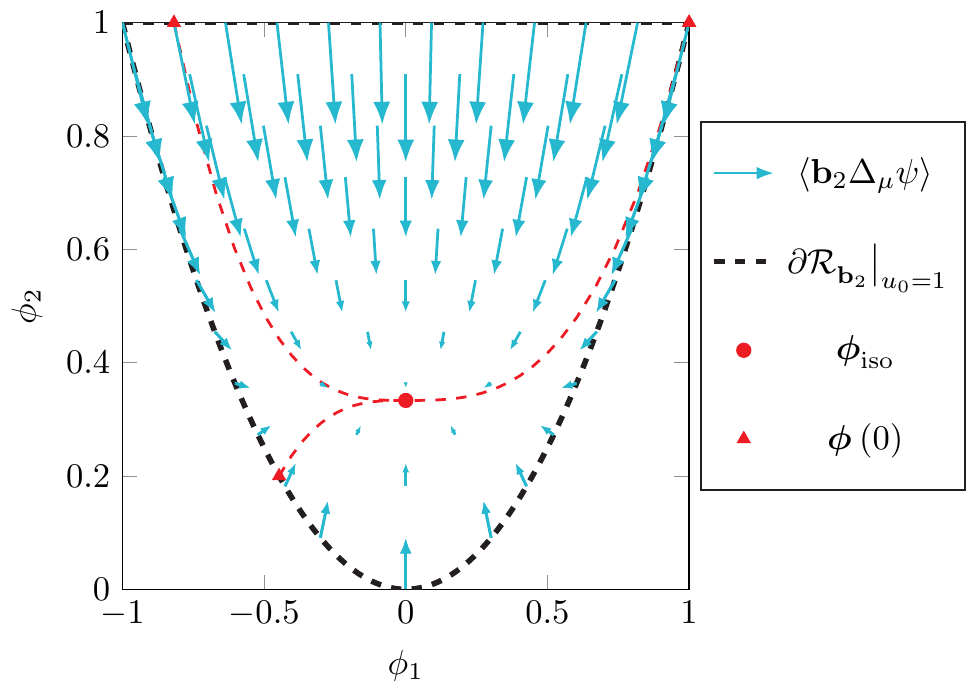}{\relinput{Images/M2LaplaceBeltramiVectorField}}
\caption{Vector field of the right-hand side and some solution trajectories of \eqref{eq:LaplaceBeltramiMomentsEquationFM} for $\momentorder = 2$ and $\momentcomp{0} = 1$. The length of the arrows is scaled by $0.03$.}
\label{fig:M2LaplaceBeltramiVectorField}
\end{figure}

Some solution curves (red dotted), starting at the realizability boundary (red triangles), are shown as well. All those curves end in the isotropic point (red dot) implying that the stationary solution of \eqref{eq:CollisionDistributionEquation} is recovered. This is not by accident. The solution of \eqref{eq:LaplaceBeltramiMomentsEquation} with a full-moment basis turns out to be just the projection of \eqref{eq:CollisionDistributionEquation} onto the corresponding moment space. This is proven below in \lemref{lem:LaplaceBeltramiFMStayRealizable}.\\

As visible in \figref{fig:M2LaplaceBeltramiVectorField}, the vector field in $\pm\normalizedmomentcomp{1} = \normalizedmomentcomp{2} = 1$ is tangential to the realizability boundary $\dRDone{\fmbasis[2]}$. Therefore, no explicit time discretization of \eqref{eq:LaplaceBeltramiMomentsEquation} generally preserves realizability using a fixed non-negative time step.

This can be shown by a simple calculation. Due to \eqref{eq:CollisionPropertyMass}, it suffices to choose $\momentcomp{0} = 1$ and therefore the explicit discretization with step size $\dtstepsize$ in normalized moments reads
\begin{align*}
\normalizedmomentcomp{1}(\timevar+\dtstepsize) &= \normalizedmomentcomp{1}(\timevar)-2\dtstepsize\normalizedmomentcomp{1}(\timevar),\\
\normalizedmomentcomp{2}(\timevar+\dtstepsize) &= \normalizedmomentcomp{2}(\timevar)-6\dtstepsize\normalizedmomentcomp{2}(\timevar)+2\dtstepsize.
\end{align*}
Plugging in $\normalizedmoments(\timevar) = \left(1,1\right)$, the updated normalized moment is given by 
\begin{align*}
\normalizedmomentcomp{1}(\timevar+\dtstepsize) &= 1-2\dtstepsize,\\
\normalizedmomentcomp{2}(\timevar+\dtstepsize) &= 1-4\dtstepsize.
\end{align*}
The update $\normalizedmoments(\timevar+\dtstepsize)$ is realizable (see \eqref{eq:M2Realizability2}) if
\begin{align*}
1 \geq \normalizedmomentcomp{2}(\timevar+\dtstepsize) = 1-4\dtstepsize\geq \normalizedmomentcomp{1}(\timevar+\dtstepsize)^2 = \left(1-2\dtstepsize\right)^2 = 1-4\dtstepsize+4\dtstepsize^2.
\end{align*}
The last inequality implies $4\dtstepsize^2\leq 0$, which is for $\dtstepsize\in\R$ only possible if $\dtstepsize=0$.
\end{example}

\begin{remark}
This is in contrast to the linear collision operator \eqref{eq:collisionOperatorR}, which is in principle easy to control since its moments are always of the form
\begin{align*}
\ints{\basis \collision{\distribution}} = \collisionrealizablepart-\moments,
\end{align*}
where $\collisionrealizablepart\in\RD{\basis}{}$ \cite{Schneider2015a}. Note that this is true for any angular basis, not only for full moments. Since the realizable set is a convex cone, this additional realizable term does not affect the realizability of the moment systems solution in a negative way, even if everything is discretized explicitly. The explicit update for \eqref{eq:LaplaceBeltramiMomentsEquation} reads
\begin{align*}
\moments(\timevar+\dtstepsize) = \left(1-\dtstepsize\right)\moments(\timevar)+\dtstepsize\collisionrealizablepart,
\end{align*}
which is realizable as long as $0\leq \dtstepsize\leq 1$. This corresponds to the standard stability condition for the explicit Euler scheme and depends on the stiffness of the system under consideration.
\end{remark}

% % % % % % % % % % % % % % % % % % % % % %
As a consequence of \assref{ass:CollisionOperator}(3), the solution $\distribution$ of \eqref{eq:CollisionDistributionEquation} is non-negative. Using this information one can conclude realizability of the exact solution of 
\begin{align}
\label{eq:ReducedMomentsEquation}
\dt \moments = \ints{\basis\collision{\ansatz[\moments]}} =: \collisionu{\moments},
\end{align}
under the following assumptions.

\begin{assumption}
\label{ass:MomentSystem}
\begin{enumerate}[(a)]
\item The map $\moments\to\ints{\basis\collision{\ansatz[\moments]}}$ is Lipschitz-continuous in $\moments$ (with respect to any norm in $\R^\momentnumber$)
\item \eqref{eq:ReducedMomentsEquation} admits a unique solution $\moments(\timevar)$ for all $\timevar\geq 0$.
\end{enumerate}
\end{assumption}

\begin{lemma}
\label{lem:LaplaceBeltramiFMStayRealizable}
\mbox{ }\\
Let $\moments(0)\in\RD{\basis}{}$ and \assref{ass:MomentSystem} be valid. Then, the solution $\moments(\timevar)$ of \eqref{eq:ReducedMomentsEquation} satisfies $\moments(\timevar)\in\RD{\basis}{}$ for all $\timevar\geq 0$.
\end{lemma}
\begin{proof}
%Since the mapping $\moments\to\collisionu{\moments}$ is Lipschitz-continuous, the solution of \eqref{eq:ReducedMomentsEquation} is unique due to the Picard-Lindel\"of theorem \cite{Wilke2011}.
Let $\distribution(\timevar,\SCheight)$ denote the solution of \eqref{eq:CollisionDistributionEquation}. As mentioned before, $\distribution(\timevar,\SCheight)\geq 0$ for all $\timevar\geq 0$ and $\SCheight\in[-1,1]$. Defining the moments of $\distribution$ as $\moments[\distribution] = \ints{\basis\distribution}$, it is immediately obvious that $\moments[\distribution]$ also solves \eqref{eq:ReducedMomentsEquation} and $\moments[\distribution](\timevar)\in\RD{\basis}{}$ for all $\timevar\geq 0$. Due to the uniqueness of the solution of \eqref{eq:ReducedMomentsEquation} (\assref{ass:MomentSystem}(b)) it follows that $\moments = \moments[\distribution]$, which completes the proof.
\end{proof}

Consequently, an implicit discretization of the moment system preserves realizability.
\begin{corollary}
\label{cor:ImplicitDiscretization}
Let $\moments(0)\in\RD{\basis}{}$. Then the implicit time-discretization 
\begin{align}
\label{eq:ImplicitDiscretization}
 \moments(\timevar+\dtstepsize) =  \moments(\timevar)+\dtstepsize \collisionu{\moments(\timevar+\dtstepsize)}
\end{align}
of \eqref{eq:ReducedMomentsEquation} satisfies $\moments(\timevar)\in\RD{\basis}{}$ for all $\timevar = j\dtstepsize$, $j\in\N$.
\end{corollary}
\begin{proof}
Similar to the proof of \lemref{lem:LaplaceBeltramiFMStayRealizable}, one can make use of the discretization of the kinetic equation \eqref{eq:CollisionDistributionEquation}, which reads 
\begin{align*}
\distribution(\timevar+\dtstepsize,\SCheight) =  \distribution(\timevar,\SCheight)+\dtstepsize \collision{\distribution(\timevar+\dtstepsize,\SCheight)}.
\end{align*}
Using \eqref{eq:CollisionDistributionPositivity} it follows that $\distribution(\timevar+\dtstepsize,\SCheight)\geq 0$, since by assumption $\distribution(\timevar,\SCheight)\geq 0$.

The solution of the system \eqref{eq:ImplicitDiscretization} is unique by Banach's fixed point theorem (using a norm that is suitably scaled by the Lipschitz constant of $\collisionopu$). As above, this solution has to satisfy $ \moments(\timevar+\dtstepsize) = \ints{\basis \distribution(\timevar+\dtstepsize)}$ and is therefore realizable.
\end{proof}

\begin{example}
We want to show that the Laplace-Beltrami operator satisfies \eqref{eq:CollisionDistributionPositivity}.

Assuming that at time $\timevar$ the solution is non-negative, the implicit discretization of \eqref{eq:CollisionDistributionEquation} can be written as
\begin{align}
\label{eq:DistributionImplicitDiscretization}
(I-\dtstepsize\LaplaceBeltramiProjection) \distribution(\timevar+\dtstepsize) =  \distribution(\timevar) \geq 0.
\end{align}
Since the Laplace-Beltrami operator is a negative operator, the operator $(I-\dtstepsize\LaplaceBeltramiProjection)$ is positive and consequently $\distribution(\timevar+\dtstepsize)\geq 0$. This can be derived rigorously by defining the Hilbert space 
\begin{align*}
\SobolevToyproblem = \{\SobolevTestfunction\in\Lp{2}(-1,1)~|~\sqrt{1-\SCheight^2}\cfrac{d\SobolevTestfunction}{d\SCheight}\in\Lp{2}(-1,1) \}
\end{align*}
with the inner product 
\begin{align*}
\left(\SobolevTestfunction,\distribution\right)_\SobolevToyproblem = \ints{\SobolevTestfunction\distribution + \dtstepsize(1-\SCheight^2)\cfrac{dv}{d\SCheight}\cfrac{d\distribution}{d\SCheight}}
\end{align*}
and the induced norm $\norm{\SobolevTestfunction}{\SobolevToyproblem} = \sqrt{\left(\SobolevTestfunction,\SobolevTestfunction\right)_\SobolevToyproblem}$. These definitions roughly follow \cite{Degond1987}. The weak formulation of \eqref{eq:DistributionImplicitDiscretization} reads
\begin{align*}
\left(\SobolevTestfunction,\distribution(\timevar+\dtstepsize)\right)_\SobolevToyproblem = \ints{\SobolevTestfunction\distribution(\timevar)}.
\end{align*}
Choosing $\SobolevTestfunction = \distribution^-(\timevar+\dtstepsize) = \min\left(0,\distribution(\timevar+\dtstepsize)\right)$, the weak formulation turns to 
\begin{align*}
\norm{\distribution^-(\timevar+\dtstepsize)}{\SobolevToyproblem}^2 = \ints{\underbrace{\distribution^-(\timevar+\dtstepsize)}_{\leq 0}\underbrace{\distribution(\timevar)}_{\geq 0}} \leq 0.
\end{align*}
Therefore, $\distribution^-(\timevar+\dtstepsize)\equiv 0$ almost everywhere and consequently $\distribution(\timevar+\dtstepsize)\geq 0$ almost everywhere.
\end{example}

\begin{remark}
We want to remark that both collision operators \eqref{eq:LaplaceBeltrami} and \eqref{eq:collisionOperatorR} with the full-moment basis satisfy all the previous assumptions since in both cases the operator $\collisionu{\moments}$ is linear in $\moments$.
\end{remark}

%% file: Sections/scheme.tex
\section{Realizability-preserving first-order scheme}
\label{sec:RPFO}
It is easy to show that a standard explicit, first-order finite-volume scheme for \eqref{eq:MomentSystemClosed} with suitably-chosen numerical fluxes automatically preserves the realizability of the underlying solution under a CFL-type constraint if the moments of the collision operator can be written as $\ints{\basis \collision{\distribution}} = \collisionrealizablepart-\moments$, where $\collisionrealizablepart\in\RD{\basis}{}$ is realizable (see e.g. \cite{Schneider2015a,Schneider2015b,Hauck2010}).\\

Unfortunately, this is in general not possible for the Laplace-Beltrami operator. As has been shown above, an explicit discretization of the right-hand-side of \eqref{eq:MomentSystemClosed} can lead to unrealizable moments even in the rather simple case of the full-moment $\MN[2]$ model. This results from the fact that the vector field defined by $\ints{\basis \collision{\distribution}}$ can point tangential to the realizability boundary and can be avoided using an implicit discretization.

On the other hand, the hyperbolic flux, which is non-linear and usually expensive to calculate, is typically non-stiff. An implicit discretization is therefore undesired.\\

To overcome this, we treat the two parts separately using an implicit-explicit time-stepping.

In the following, the spatial domain $\Domain = (\zL, \zR)$ is divided into (for notational simplicity) $\ncells$ (equidistant) cells $\cell{\cellind} = (\z_{\cellind-\frac12}, \z_{\cellind+\frac12})$, where the cell interfaces are given by $\z_{\cellind\pm\frac12} = \z_\cellind \pm 
\frac{\dzstepsize}{2}$ for cell centres $\z_\cellind = \zL + (\cellind - \frac12)\dzstepsize$, and
$\dzstepsize = \frac{\zR - \zL}{\ncells}$. 

Defining the averaging operator
\begin{align*}
\cellmean{\, \cdot \,} := \frac{1}{\dzstepsize} \int_{\cell{\cellind}} \cdot~d\z,
\end{align*}
the discretized form of \eqref{eq:GeneralHyperbolicSystem} reads
\begin{align}
\label{eq:discretizedform}
\cfrac{\momentscellmeantime{\cellind}{\timeind+1}-\momentscellmeantime{\cellind}{\timeind}}{\dtstepsize} = - \frac{1}{\dzstepsize}\left(\numericalFlux(\momentscellmeantime{\cellind}{\timeind}, \momentscellmeantime{\cellind+1}{\timeind})-\numericalFlux(\momentscellmeantime{\cellind-1}{\timeind}, \momentscellmeantime{\cellind}{\timeind})\right)
+ \Source\left(\momentscellmeantime{\cellind}{\timeind+1}\right)
\end{align}
where $\numericalFlux$ is a numerical flux function coupling the solution on cell $\cell{\cellind}$ with its neighbours.

We use the kinetic flux (see e.g. \cite{Schneider2015b,Hauck2010,Garrett2014})
\begin{align}
\label{eq:GodunovFlux1}
 \numericalFlux(\moments[1], \moments[2]) &= \ints{\basis \monotoneFlux{\distributionpv{1}}{\distributionpv{2}}}, \qquad \text{ where}\\
\monotoneFlux{\distributionpv{1}}{\distributionpv{2}} &= \begin{cases}
\SCheight\distributionpv{1} & \text{ if } \SCheight\geq 0,\\
\SCheight\distributionpv{2} & \text{ if } \SCheight\leq 0
\end{cases}
\end{align}
and $\distributionpv{1,2}$ are the ans\"atze for $\moments[1,2]$, respectively. This is generally possible for minimum-entropy and Kershaw models by carrying out the integrations over the half-spaces separated by $\SCheight = 0$. This is particularly easy in case of half- and mixed-moment models since then the numerical flux can be explicitly written in terms of the moments instead of some half moments of the ansatz function.

The incorporation of boundary conditions is non-trivial. Here, an often-used approach is taken that incorporates boundary conditions via `ghost cells'.
First assume that it is possible to smoothly extend $\distributionboundary(\timevar,\z,\SCheight)$ in $\SCheight$ to $[-1,1]$ for $\z\in\{\zL,\zR\}$ (note that while moments are defined using integrals over all $\SCheight$, the
boundary conditions in \eqref{eq:TransportEquation1DBCa}--\eqref{eq:TransportEquation1DBCb} are only defined for
$\SCheight$ corresponding to incoming data). 

Then the moment approximations in the ghost cells at $\z_0$ and $\z_{\ncells+1}$ simply take the form 
\begin{subequations}
\label{eq:BCMomentSystemSimple1D}
\begin{align}
\momentslocal{0}(\timevar, \z_{\frac12}) &:= \ints{\basis \distributionboundary(\timevar,\zL,\SCheight)},\\
\momentslocal{\ncells + 1}(\timevar, \z_{\ncells + \frac12}) &:= \ints{\basis \distributionboundary(\timevar,\zR,\SCheight)}.
\end{align}
\end{subequations}
Note, however, that the validity of this approach, due to its inconsistency with 
the original boundary conditions  \eqref{eq:TransportEquation1DBCa}--\eqref{eq:TransportEquation1DBCb}, is not entirely non-controversial, but the question of appropriate boundary conditions for
moment models is an open problem \cite{pomraning1964variational,Larsen1991,Rulko1991,Struchtrup2000,levermore2009boundary} which is not explored here. 

The IMEX time-stepping in \eqref{eq:discretizedform} uses the forward-backward Euler scheme \cite{Ascher1997}. Since this is nothing else than doing a Godunov splitting of the hyperbolic part (treated explicitly) and the (stiff) source term (treated implicitly), the following theorem can be concluded.

\begin{theorem}
\label{thm:RP}
Let $\momentscellmeantime{\cellind}{\timeind}\in\RD{\basis}{}$ for all $\cellind=0,\ldots,\ncells+1$. Furthermore, let \assref{ass:CollisionOperator} and \assref{ass:MomentSystem} hold, and $\absorption(\timevar,\x),\scattering(\timevar,\x),\source(\timevar,\x,\SCheight)\in\Rpos$ be bounded and continuous in $\timevar$. 

Then, the IMEX scheme \eqref{eq:discretizedform} preserves realizability (i.e. $\momentscellmeantime{\cellind}{\timeind+1}\in\RD{\basis}{}$ for all $\cellind=1,\ldots,\ncells$) under the CFL condition 
\begin{align}
\label{eq:CFL}
 \dtstepsize \leq \dzstepsize.
\end{align}
\end{theorem}
\begin{proof}
The scheme \eqref{eq:discretizedform} is equivalent to the following splitting scheme
\begin{subequations}
\begin{align}
\label{eq:discretizedform2a}
\momentscellmeantime{\cellind}{*} &= \momentscellmeantime{\cellind}{\timeind}- \frac{\dtstepsize}{\dzstepsize}\left(\numericalFlux(\momentscellmeantime{\cellind}{\timeind}, \momentscellmeantime{\cellind+1}{\timeind})-\numericalFlux(\momentscellmeantime{\cellind-1}{\timeind}, \momentscellmeantime{\cellind}{\timeind})\right),\\
\label{eq:discretizedform2b}
\momentscellmeantime{\cellind}{\timeind+1} &= \momentscellmeantime{\cellind}{*}+\dtstepsize\Source\left(\momentscellmeantime{\cellind}{\timeind+1}\right).
\end{align}
\end{subequations}

We recapitulate the arguments from e.g. \cite{Schneider2015a,Schneider2015b} to show that \eqref{eq:discretizedform2a} preserves realizability. We have that 
\begin{align*}
\momentscellmeantime{\cellind}{*} &= \ints{\distribution[*]}\\
\distribution[*] &= \ansatz[\momentscellmeantime{\cellind}{\timeind}]-\frac{\dtstepsize}{\dzstepsize}\left(\max\left(\SCheight,0\right)\left(\ansatz[\momentscellmeantime{\cellind}{\timeind}]-\ansatz[\momentscellmeantime{\cellind-1}{\timeind}]\right)+\min\left(\SCheight,0\right)\left(\ansatz[\momentscellmeantime{\cellind+1}{\timeind}]-\ansatz[\momentscellmeantime{\cellind}{\timeind}]\right)\right)\\
&\geq \left(1-\frac{\dtstepsize}{\dzstepsize}\right)\ansatz[\momentscellmeantime{\cellind}{\timeind}]\stackrel{\eqref{eq:CFL}}{\geq} 0,
\end{align*}
where $\ansatz[\momentscellmeantime{\cellind-1}{\timeind}],\ansatz[\momentscellmeantime{\cellind}{\timeind}],\ansatz[\momentscellmeantime{\cellind+1}{\timeind}]\geq 0$ are the respective ans\"atze \eqref{eq:psiME} for the moment vectors in cells $\cell{\cellind-1},\cell{\cellind}$ and $\cell{\cellind+1}$. Thus, $\momentscellmeantime{\cellind}{*}$ is generated by the non-negative distribution function $\distribution[*]$ and is therefore realizable, i.e. $\momentscellmeantime{\cellind}{*}\in\RD{\basis}{}$.

To show a similar result for \eqref{eq:discretizedform2b}, \corref{cor:ImplicitDiscretization} has to be adopted to the situation. The update has the form 
\begin{align*}
\momentscellmeantime{\cellind}{\timeind+1} &= \left(\momentscellmeantime{\cellind}{*}+\dtstepsize\ints{\basis\cellmean[\cellind]{\source}}\right)+\dtstepsize\left(\cellmean[\cellind]{\scattering}\collisionu{\momentscellmeantime{\cellind}{\timeind+1}}-\cellmean[\cellind]{\absorption}\momentscellmeantime{\cellind}{\timeind+1}\right)\\
&=\underbrace{\ints{\basis\underbrace{\left(\ansatz[\momentscellmeantime{\cellind}{*}]+\dtstepsize\cellmean[\cellind]{\source}\right)}_{\geq 0}}}_{\in\RD{\basis}{}}+\dtstepsize\left(\cellmean[\cellind]{\scattering}\collisionu{\momentscellmeantime{\cellind}{\timeind+1}}-\cellmean[\cellind]{\absorption}\momentscellmeantime{\cellind}{\timeind+1}\right).
\end{align*}
This can be stated equivalently as 
\begin{align*}
\left(1+\dtstepsize\cellmean[\cellind]{\absorption}\right)\momentscellmeantime{\cellind}{\timeind+1}
&=\ints{\basis\left(\ansatz[\momentscellmeantime{\cellind}{*}]+\dtstepsize\cellmean[\cellind]{\source}\right)}+\dtstepsize\cellmean[\cellind]{\scattering}\collisionu{\momentscellmeantime{\cellind}{\timeind+1}}.
\end{align*}
Since the realizable set is a convex cone and $\cellmean[\cellind]{\absorption}\geq 0$, $\left(1+\dtstepsize\cellmean[\cellind]{\absorption}\right)^{-1}\ints{\basis\left(\ansatz[\momentscellmeantime{\cellind}{*}]+\dtstepsize\cellmean[\cellind]{\source}\right)}\in\RD{\basis}{}$. Thus, 
\begin{align*}
\momentscellmeantime{\cellind}{\timeind+1}
&=\left(1+\dtstepsize\cellmean[\cellind]{\absorption}\right)^{-1}\ints{\basis\left(\ansatz[\momentscellmeantime{\cellind}{*}]+\dtstepsize\cellmean[\cellind]{\source}\right)}+\left(1+\dtstepsize\cellmean[\cellind]{\absorption}\right)^{-1}\dtstepsize\cellmean[\cellind]{\scattering}\collisionu{\momentscellmeantime{\cellind}{\timeind+1}}
\end{align*}
is of the form that \corref{cor:ImplicitDiscretization} can be applied (with suitable redefinitions of $\dtstepsize$). Note that boundedness and continuity of the physical parameters are necessary such that a similar modification of \lemref{lem:LaplaceBeltramiFMStayRealizable} is still valid.

Thus $\momentscellmeantime{\cellind}{\timeind+1}\in\RD{\basis}{}$, which completes the proof.
\end{proof}

\begin{remark}
Using an explicit discretization of the source term, the CFL condition \eqref{eq:CFL} has to be modified to 
\begin{align*}
\dtstepsize\leq \frac{1}{\frac{1}{\dzstepsize}+\max\limits_{\cellind,\timeind}\left(\cellmean[\cellind]{\absorption}\left(\timevar_\timeind\right)+\cellmean[\cellind]{\scattering}\left(\timevar_\timeind\right)\right) }
\end{align*}
to preserve realizability (if possible at all) \cite{Schneider2015a,Schneider2015b}.
\end{remark}

%% file: Sections/results.tex
\section{Numerical experiments}
\label{sec:NumExp}
\subsection{Manufactured solution}
\label{sec:manu-soln}
In general, analytical solutions for minimum-entropy models are not known.
Therefore, to test the convergence and efficiency of our scheme, the
method of manufactured solutions is used, following the target solution given
in \cite{Schneider2015b}. 
The solution is defined on the spatial domain $\Domain = (-\pi, \pi)$ with
periodic boundary conditions.

A kinetic density in the form of the entropy ansatz is given by
\begin{align}
 \analyticalsolution(\timevar,\z,\SCheight) =& \exp\left(\multiplierscomp{0}(\timevar,\z) + \multiplierscomp{1}(\timevar,\z) \SCheight \right), \label{eq:MFSM3}\\
 \multiplierscomp{0}(\timevar,\z) =& -\MFSconstK - \sin(\z-\timevar) - \MFSconsta,\nonumber\\
 \multiplierscomp{1}(\timevar,\z) =&  \MFSconstK + \sin(\z-\timevar).\nonumber
\end{align}
A source term is defined by applying the transport operator to $\analyticalsolution$, giving
$$
\source(\timevar,\z,\SCheight) := \dt \analyticalsolution(\timevar,\z,\SCheight) + \SCheight \dz
\analyticalsolution(\timevar,\z,\SCheight) + \absorption(\timevar, \z) \analyticalsolution(\timevar,\z,\SCheight),
$$
where
$$
 \absorption(\timevar, \z) := 4\left(1 - \cos\left(\z - \timevar\right)\right).
$$
Thus, by inserting this $\source$ into \eqref{eq:TransportEquation1D} and setting 
$\scattering = 0$,
$\analyticalsolution$ is a solution of
\eqref{eq:TransportEquation1D}.

A straightforward computation shows that $\source \ge 0$ (for
any $\MFSconsta$ and $\MFSconstK$), which means that \thmref{thm:RP} can be applied to the
resulting moment system.

Furthermore, $\MFSconsta$ is chosen as
\begin{align*}
 \MFSconsta &= - \MFSconstK + 1 - \log\left( \cfrac{\MFSconstK - 1}{2\sinh(\MFSconstK - 1)} \right)
\end{align*}
so that the maximum value of $\ints{\analyticalsolution}$ for $(\timevar, \z) \in [0, \tf] \times \Domain$
is one.
As $\MFSconstK$ is increased, $\analyticalsolution$ converges to a Dirac delta at $\SCheight = 1$.\\

Since $\analyticalsolution$ has the form of an entropy ansatz,
$\analyticalmoments = \ints{\basis \analyticalsolution}$ is also a solution of \eqref{eq:MomentSystemUnclosed1D}
whenever $1$ and $\SCheight$ are in the linear span of the basis $\basis$. Notice also that $\analyticalmoments$ approaches the boundary of realizability as $\MFSconstK$ is
increased.\\

The final time is chosen to be $\tf = \pi / 5$ while $\MFSconstK\in\{2,25\}$ is used, for which the normalized first-order moment satisfies $\frac{\analyticalmomentcomp{1}}{\analyticalmomentcomp{0}} \in \{[0.313,0.672],[0.958,0.962]\}$
(recall that $\abs{\analyticalmomentcomp{1}} \leq \analyticalmomentcomp{0}$ is necessary for realizability). 

In the following, the $\MN[3]$ model is used so that the results include the
effects of the numerical optimization.

Errors are computed in the zeroth moment of the solution $\analyticalmomentcomp{0}(\timevar, \z) := \ints{\analyticalsolution(\timevar, \z, \cdot)}$.
Then $\Lp{1}$- and $\Lp{\infty}$-errors for the zeroth moment $\momentcompprojected{0}(\timevar, \z)$
(that is, the zeroth component of a numerical solution $\momentsprojected$) are
defined as
\begin{equation}
 \LpError{1} = \dzstepsize\sum\limits_{\cellind=1}^{\ncells} \left|\cellmean[\cellind]{\analyticalmomentcomp{0}}(\tf) - \cellmean[\cellind]{\momentcomp{0}}(\tf) \right|
  \quand
\LpError{\infty} = \max_{\cellind=1,\ldots,\ncells} \left|\cellmean[\cellind]{\analyticalmomentcomp{0}}(\tf) - \cellmean[\cellind]{\momentcomp{0}}(\tf) \right|,
\label{eq:errors}
\end{equation}
respectively. The observed convergence order $\MFSconstorder$ is defined by
\begin{equation}
 \frac{\LpError[h_1]{\MFSp}}{\LpError[h_2]{\MFSp}} = \left( \frac{\dzstepsize_1}{\dzstepsize_2} \right)^\MFSconstorder,
\label{eq:conv-order}
\end{equation}
where $\LpError[h_i]{\MFSp}$, $i \in \{1, 2\}$, $\MFSp \in \{1, \infty\}$, is the $\Lp{\MFSp}$-error $\LpError{\MFSp}$ for the
numerical solution using cell size $\dzstepsize_i$.

A convergence table for two different values of $\MFSconstK$ is presented in \tabref{tab:ConvergenceDG}. They correspond in spatial average to the sets of normalized moments $\normalizedmoments = \left(0.515,0.463,0.333\right)^T$ ($\MFSconstK = 2$) and $\normalizedmoments = \left(0.960,0.923,0.889\right)^T$ ($\MFSconstK = 25$) with relative distance to the realizability boundary (absolute distance divided by the maximal possible distance) of $5.016\%$ and $0.0006\%$, respectively.

\begin{table}[htbp]
\centering
\begin{tabular}{r r@{.}l c r@{.}l c r@{.}l c r@{.}l c }
& \multicolumn{6}{c}{$\MFSconstK = 2 $}& \multicolumn{6}{c}{$\MFSconstK = 25 $}\\
\cmidrule(r){2-7} \cmidrule(r){8-13} 
$\ncells $ & \multicolumn{2}{c}{$E^1_h$} & $\nu$ & \multicolumn{2}{c}{$E^\infty_h$} & $\nu$& \multicolumn{2}{c}{$E^1_h$} & $\nu$ & \multicolumn{2}{c}{$E^\infty_h$} & $\nu$\\ %\midrule 
\cmidrule(r){1-1} \cmidrule(r){2-4} \cmidrule(r){5-7} \cmidrule(r){8-10} \cmidrule(r){11-13}
 40 & 5 & 332e-02 & ---& 2 & 355e-02 & ---& 7 & 063e-03 & ---& 2 & 633e-03 & ---\\
80 & 2 & 713e-02 & 0.97& 1 & 208e-02 & 0.96& 3 & 558e-03 & 0.99& 1 & 329e-03 & 0.99\\
160 & 1 & 368e-02 & 0.99& 6 & 118e-03 & 0.98& 1 & 792e-03 & 0.99& 6 & 671e-04 & 0.99\\
320 & 6 & 862e-03 & 1.00& 3 & 078e-03 & 0.99& 9 & 035e-04 & 0.99& 3 & 341e-04 & 1.00\\
640 & 3 & 444e-03 & 0.99& 1 & 554e-03 & 0.99& 4 & 655e-04 & 0.96& 1 & 684e-04 & 0.99\\
%1280 & 1 & 766e-03 & 0.96& 7 & 810e-04 & 0.99\\
\end{tabular}
\caption{$\Lp{1}$- and $\Lp{\infty}$-errors and observed convergence order $\MFSconstorder$
for the IMEX kinetic scheme with $\MN[3]$ manufactured solution \eqref{eq:MFSM3} and optimization
gradient tolerance $\opttol = 10^{-6}$.}
\label{tab:ConvergenceDG}
\end{table}

It can be observed that the expected convergence rates are achieved both in $\Lp{1}$-
and $\Lp{\infty}$-errors.

\begin{remark}
The scheme is not convergent for arbitrarily large values of $\MFSconstK$.
For big $\MFSconstK$, the numerical solution veers so close to the
boundary of the realizable set that the optimization has to use
regularization, thus introducing errors into the solution. This has been shown in \cite{Schneider2015a} for a simpler convergence test and was also observed before in \cite{Alldredge2014}.
\end{remark}

%We use the one-dimensional version of the first-order, realizability-preserving, implicit-explicit kinetic scheme derived in \cite{Schneider2015c}. All results are computed on a grid with $1000$ points. The reference solution is given by the $\PN[99]$ model.
%
\subsection{Plane source}
\label{sec:Planesource}
In this test case an isotropic distribution with all mass concentrated in the middle of an infinite domain $\z \in
(-\infty, \infty)$ is defined as initial condition, i.e.
\begin{align*}
 \distributiontzero(\z, \SCheight) = \distributionvacuum + \delta(\z),
\end{align*}
where the small parameter $\distributionvacuum = 0.5 \times 10^{-8}$ is used to
approximate a vacuum.
In practice, a bounded domain must be used which is large
enough that the boundary should have only negligible effects on the
solution. For the final time $\tf = 1$, the domain is set to $\Domain = [-1.2, 1.2]$ (recall that for all presented models the maximal speed of propagation is bounded in absolute value by one \cite{Schneider2015a,Levermore1998}).

At the boundary the vacuum approximation
\begin{align*}
 \distributionboundary(\timevar,\zL,\SCheight) \equiv \distributionvacuum \quand
 \distributionboundary(\timevar,\zR,\SCheight) \equiv \distributionvacuum
\end{align*}
 is used again. Furthermore, the physical coefficients are set to $\scattering \equiv 1$, $\absorption \equiv 0$ and $\source \equiv 0$.

\begin{figure}[htbp]
 \centering
 \settikzlabel{fig:Planesource.M1D}
  \settikzlabel{fig:Planesource.M1FM}
 \settikzlabel{fig:Planesource.M2D}
  \settikzlabel{fig:Planesource.M2FM}
 \settikzlabel{fig:Planesource.M3D}
  \settikzlabel{fig:Planesource.M3FM}    
\externaltikz{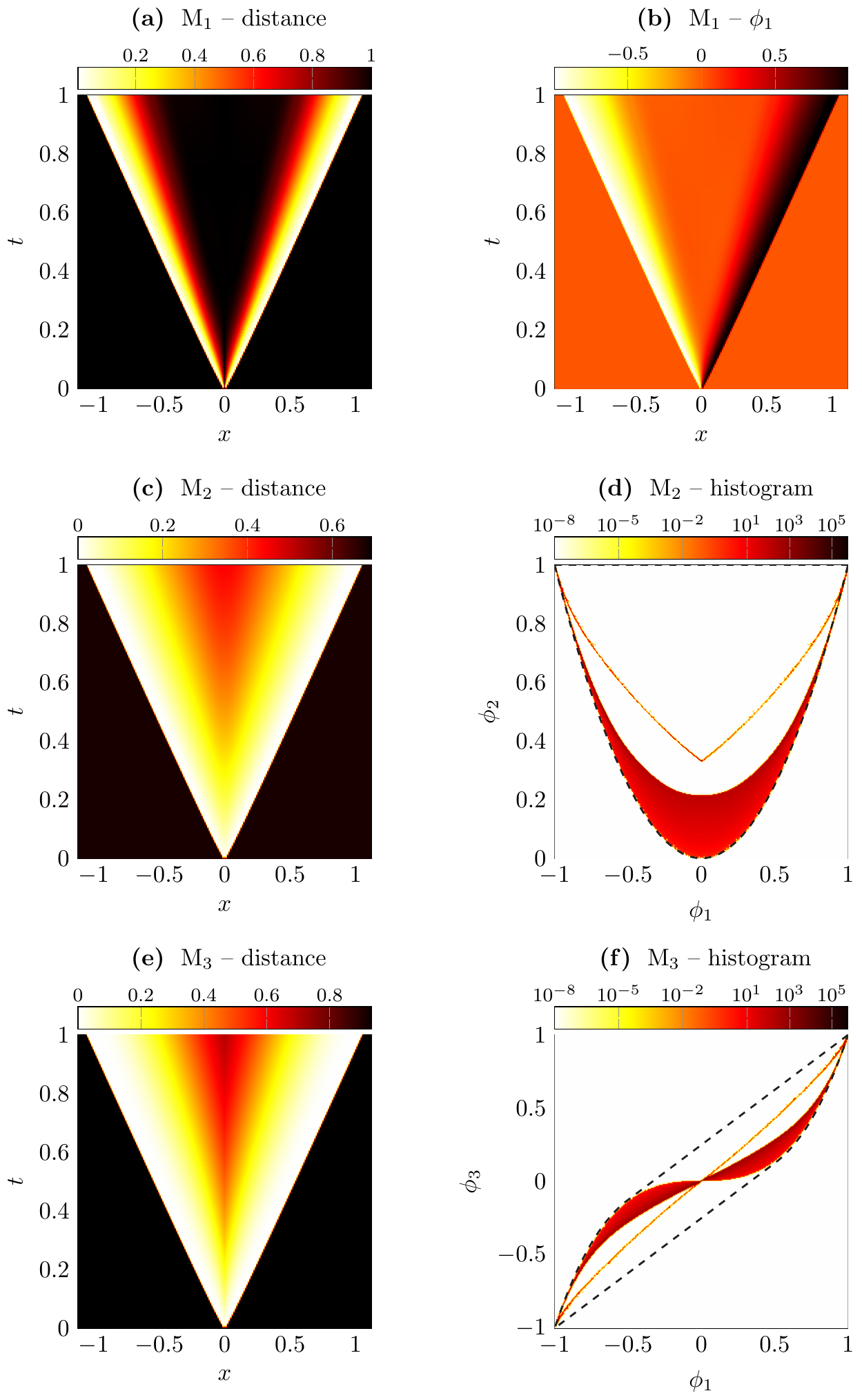}{\relinput{Images/PlanesourceLB}}
\caption{Relative distance to the realizability boundary and related quantities for the plane-source test.}
 \label{fig:Planesource}
\end{figure}

All solutions are computed with an even number of cells, so the initial Dirac delta lies on a cell boundary.
Therefore it is approximated by splitting it into the cells immediately to the left and right. In \figref{fig:Planesource}, only positive $\z$ are shown since the solutions are always symmetric around $\z = 0$.\\

\begin{figure}[htbp]
\externaltikz{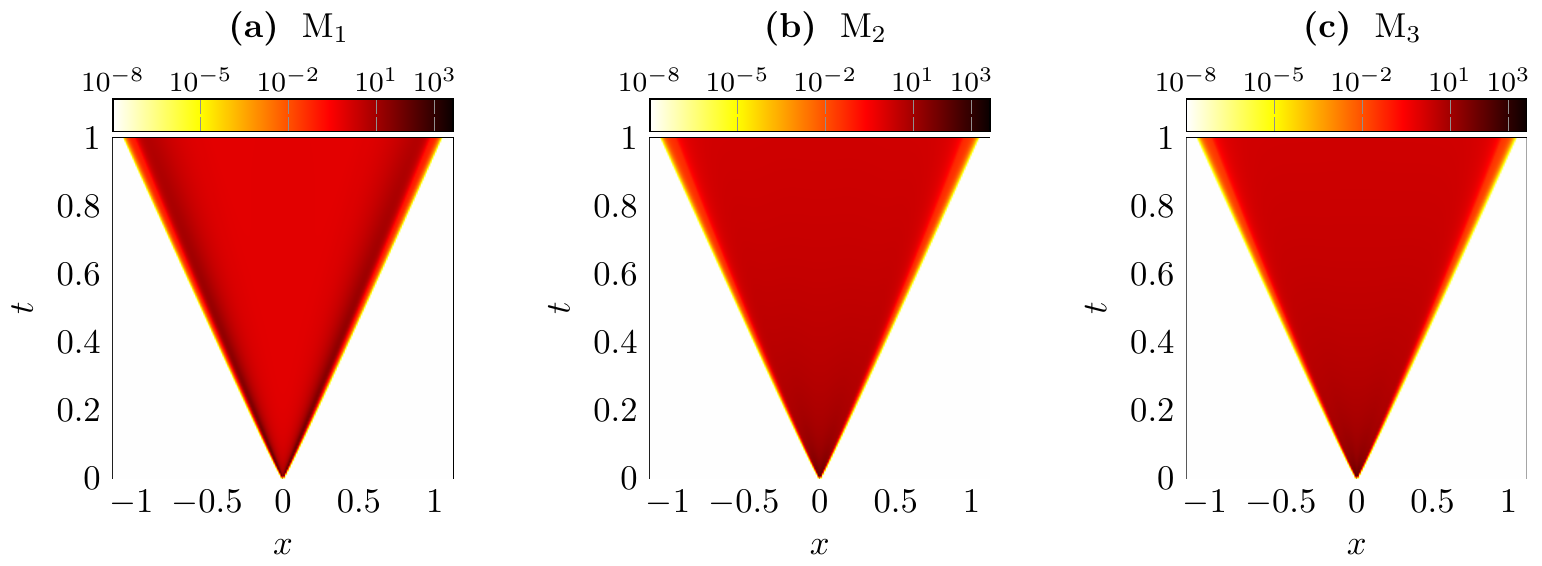}{\relinput{Images/PlanesourceCuts}}
 \centering
\caption{Results for the plane-source test in the space-time domain in a logarithmic scale.}
 \label{fig:PlanesourceCuts}
\end{figure}

For an intense discussion of the solution of the moment models see e.g. \cite{Schneider2014,Schneider2016}. For convenience, the space-time behaviour of the density $\density$ for $\MN[1]$ to $\MN[3]$ are shown in \figref{fig:PlanesourceCuts}.
%The figure shows the solution of the $\DMMN[2]$ model in comparison to some mixed-moment $\MMN$ and full-moment $\MN$ models with a similar number of degrees of freedom (i.e. the number of moments $\momentnumber$). 

%Observe that the difference between $\MMN[2]$ and $\DMMN[2]$ is negligible\footnote{This is no longer true if the isotropic scattering operator $\collision{\distribution} = \distribution - \frac12\int\limits_{-1}^1\distribution(\SCheight')~d\SCheight'$ is used.}. Although the $\DMMN[2]$ model is exactly between $\MMN[1]$ and $\MMN[2]$ (regarding degrees of freedom), its solution is much closer to those of the $\MMN[2]$ model.

%Doing the same comparison with the $\MN$ models shows that the $\DMMN[2]$ model is closer to the $\MN[2]$ than to the $\MN[3]$ model (while all three models differ insignificantly from the reference solution\footnote{This results from the quadratic dependence of the Laplace-Beltrami eigenvalues with respect to the moment order $\momentorder$.}).

A known problem of the minimum-entropy approach is the fact that close to the realizability boundary the moment system becomes ill-conditioned \cite{AllHau12}. We investigate the relative distance of the plane-source $\MN$ solutions to the boundary of the \emph{normalized realizable set} 
\begin{align*}
\RDone{\basis} = \left\{\moments~:~\exists \distribution(\SC)\ge 0,\, \momentcomp{0} = \ints{\distribution} =1,
 \text{ such that } \moments =\ints{\basis\distribution} \right\}
\end{align*}
as the ratio of the euclidean distance to the realizability boundary and the maximal possible distance, i.e.
\begin{align*}
\distRDrel{\moments} = \cfrac{\distRD{\moments}}{\max\limits_{\hat{\moments}\in\RDone{\basis}} \distRD{\hat{\moments}}}, \qquad \distRD{\moments} = \min\limits_{\hat{\moments}\in\dRDone{\basis}} \norm{\moments-\hat{\moments}}{2}.
\end{align*}
The maximal distances are
\begin{align*}
\max\limits_{\hat{\moments}\in\RDone{\basis}} \distRD{\hat{\moments}}&=1 && \text{for }\MN[1],\\
\max\limits_{\hat{\moments}\in\RDone{\basis}} \distRD{\hat{\moments}}&=\frac12 && \text{for }\MN[2],\text{ and}\\
\max\limits_{\hat{\moments}\in\RDone{\basis}} \distRD{\hat{\moments}}&=\frac15 && \text{for }\MN[3].
\end{align*}
The results are shown in \figref{fig:Planesource}.

While the relative distance in case of the $\MN[1]$ model (Figures~\ref{fig:Planesource.M1D} and \ref{fig:Planesource.M1FM}) is directly related the normalized first moment $\normalizedmomentcomp{1}$ by $\distRDrel{\moments} = 1-\abs{\normalizedmomentcomp{1}}$ (resulting in small distances only close to the peak at $\timevar = \pm \x$), the distances in case of the higher-order models become smaller even in the interior of the set $\{\abs{\x}\leq \timevar\}$. The minimal values that occurred are $0.0039$ ($\MN[1]$), $2.2147\cdot10^{-5}$ ($\MN[2]$) and $9.0981\cdot10^{-7}$ ($\MN[3]$), showing that the underlying moment problem becomes harder to solve with increasing moment order $\momentorder$.

Figures~ \ref{fig:Planesource.M2FM} and \ref{fig:Planesource.M3FM} show a histogram ($300\cdot300$ bins) of the $\MN[2]$ and $\MN[3]$ solution in the $\normalizedmomentcomp{1}-\normalizedmomentcomp{\momentorder}$ phase space (where $\momentorder$ is either $2$ or $3$, respectively). The histogram is built out of the solution values at the $10000$ cell centres and $100$ time frames. The boundary of the (projected) normalized realizable set is depicted as a black, dashed line. In case of the $\MN[2]$ model, it is visible that mostly the lower part of the realizable set is filled with particles, complemented with a stream of particles connected to the isotropic point $\normalizedisotropicmoment = \left(0,\frac13\right)$. No particles occur close to the point of maximal distance $\normalizedmoments=\left(0,\frac12\right)$. Thus, the relative distance for $\MN[2]$ is always strictly smaller than $1$. A similar effect occurs in case of the $\MN[3]$ model, but less pronounced.

%% file: Sections/outlook.tex
\section{Conclusions and outlook}
\label{sec:Conclusions}
We derived an implicit-explicit scheme for moment systems that are generated by a non-negative ansatz. This scheme preserves realizability under a standard CFL condition, even in the case of stiff source terms (e.g. strong scattering or absorption), while the implicit systems have to be solved only locally in every space-time cell. In many cases, these implicit systems are linear, resulting in a very efficient algorithm.

Convergence of the algorithm was tested against a manufactured solution, showing the designed first order. Furthermore, the plane-source problem served as benchmark test, showing how close to the realizability boundary the scheme can get.

While this first-order scheme is easy to implement, the benefit of high-order schemes in terms of efficiency is necessary to obtain reasonable approximations in higher dimensions in an appropriate time. Future work will investigate how to couple higher-order IMEX schemes with the fully-explicit, high-order kinetic \cite{Schneider2015b} and discontinuous-Galerkin scheme \cite{Schneider2015a}, removing the stiffness from these two methods.

Furthermore, it is unclear if the mixed-moment model \cite{Frank07,Schneider2014} in combination with the Laplace-Beltrami operator $\LaplaceBeltramiProjection$ fulfils the assumptions of \thmref{thm:RP} (it contains the microscopic quantity $\ansatz[\moments]\left(0\right)$, i.e. the solution of \eqref{eq:ReducedMomentsEquation} depends on the chosen ansatz). Nevertheless, \eqref{eq:discretizedform} performs well in practice even in this situation, which could mean that either the above assumptions are fulfilled or \thmref{thm:RP} can be extended to a weaker set of assumptions.